\documentclass[10pt]{amsart} 
\usepackage{amscd, amssymb}
\input xy
\xyoption{all}

\theoremstyle{plain}
\newtheorem{proposition}{Proposition}
\newtheorem{lemma}{Lemma}
\newtheorem{theorem}{Theorem}

\newtheorem{corollary}{Corollary}

\theoremstyle{definition}

\newtheorem*{acknowledgements}{Acknowledgements}
\newtheorem*{update}{Update}  

\theoremstyle{remark}
\newtheorem{remark}{Remark}

\newtheorem*{claim}{Claim}

\newcommand{\pp}{\mathbb{P}}

\newcommand{\qq}{\mathbb{Q}}
\newcommand{\cc}{\mathbb{C}}

\newcommand{\mbar}{{\overline{M}}}
\newcommand{\gbar}{\overline{G}}
\newcommand{\so}{{\mathcal{O}}}

\newcommand{\p}{\partial}

\newcommand{\ev}{{\operatorname{ev}}}

\newcommand{\vir}{\operatorname{vir}}

\newcommand{\rk}{\operatorname{rk}}

\newcommand{\e}{c_{\operatorname{top}}}
\newcommand{\la}{\langle}
\newcommand{\ra}{\rangle}

\newcommand{\bel}{\begin{equation} \label}
\newcommand{\ee}{{\end{equation}}}
\newcommand{\pd}{\operatorname{PD}}
\newcommand{\E}{\mathcal{E}}
\newcommand{\F}{\mathcal{F}}

\begin{document}
\title[Quantum Lefschetz Hyperplane Theorem]
      {\Large Quantum Lefschetz Hyperplane Theorem}
\author{Y.-P. Lee}
\address{Department of Mathematics\\
	UCLA\\
	Los Angeles, CA 90095-1555}
\email{yplee@math.ucla.edu}
\thanks{Research partially supported by NSF grant}


\begin{abstract}
The mirror theorem is generalized to any smooth projective variety $X$.
That is, a fundamental relation between the Gromov--Witten invariants of $X$ 
and Gromov--Witten invariants of complete intersections $Y$ in $X$ is 
established.
\end{abstract}

\maketitle


\emph{Notation.} 

\begin{list}{$\bullet$}{}
\item Our notations follow those in \cite{FP} unless otherwise mentioned.      
\item The (co)homology/Chow groups are over $\qq$ and all varieties 
(schemes, stacks) are over ground field $\cc$.
\item $X$ is a smooth projective variety. $i: X \hookrightarrow 
P:=\prod_{i=1}^N \pp^{r_i}$ is
and embedding such that $i^*: NS(\prod \pp^{r_i}) \to NS(X)$ is an
isomorphism, where $NS(X)=Z^1(X)/\sim_{\text{hom}}$ is the Neron-Severi group.
\item Let $p_i := i^* (h_i)$, where $h_i$ is the hyperplane
class of $\pp^{r_i}$. 
\item
For each $\beta \in H_2 (X)$ $d(\beta):=i_*(\beta) \in H_2(P)$. 
\item
$\pi_k: \mbar_{0,n+1}(Y, \beta) \to \mbar_{0,n}(Y,\beta)$ is the forgetful 
morphism which forgets the $k$-th marking. $\pi:= \pi_{n+1}$.
\item
$\ev: \mbar_{0,n}(X,\beta) \to (X)^n$ is the evaluation morphism.
\item
$p_k: Y_1 \times_Z Y_2 \to Y_k$ are the projection morphisms.
\item
$\gbar_{0,n}(X,\beta) := \mbar_{0,n} (X \times \pp^1, (\beta,1))$ is called 
$n$-pointed \emph{graph space} of $X$.
\item
For a convex vector bundle $E$ on $X$, define a vector bundle
$E_{\beta} := \pi^* (R^0 \pi_* \ev^* E)$ on $\mbar_{0,1}(X,\beta)$. 
If $E$ is concave, $E_{\beta} :=\pi^* (R^1 \pi_* \ev^* E)$. Similar
definition on $\gbar_{0,n}(X,\beta)$: $E^G_{\beta}:= R^j \pi_* \ev^* E$,
$j=0$ for convex $E$ and $j=1$ for concave $E$.
\item
A stable map is denoted by a quadruple $(C,x,f:C \to X,\beta)$, where $x=
(x_1, \cdots,x_n)$ are the marked points, $\beta$ is the degree of the map $f$.
\item
Let $f: V \to M$ be a morphism with $M$ smooth, then 
$f^{\vir}_* : H^*(V) \to H^*(M)$
is defined to be 
\[
 f^{\vir}_*(\alpha):= \pd \circ f_* (\alpha \cap [V]^{\vir}),
\] 
where $\pd: H_*(M) \to H^*(M)$ is the Poincar\'e duality map
and $f_* : H_*(V) \to H_*(M)$ is the push-forward of cycles. Note that
$f_*$ is also used for push-forward of cohomology classes when $V$ is
also smooth. There should be no confusion.
\item
$\psi:= c_1(\mathcal{L}_1)$ is the universal cotangent class. Namely
$\mathcal{L}=x_1^*(\omega_{C/M})$, where $\omega_{C/M}$ is the relative
dualizing sheaf of the universal curve $C \to \mbar_{0,1}(X,\beta)$ and
$x_1$ is the marked point. 
\item
All objects considered on graph space and on linear sigma model 
$\pp^r_d:=\pp^{(r+1)d+r}$ are $\cc^*$-equivariant.
\item
We have summarized the most useful notations (other than those listed above)
in the commutative diagrams \eqref{e:6} and \eqref{e:7}.
\end{list}

\parskip=1pt

\section{Introduction}
\subsection{Statement of the main results}
Let $X$ be a smooth projective variety. Let ${T_k}$ be a basis of
$H^*(X)$, $T_0=1, T_i=p_i, i=1,\cdots N$ and let $t_k$ be the dual
coordinates of $T_k$. Let $g_{mk}:=\int_X T_m \cup T_k$ the intersection
pairing and $(g^{mk})$ is the inverse matrix of $(g_{mk})$. Let $E:=
(\oplus_j L_{j})$ be a vector bundle on $X$, where $L_j$ are pull-backs of 
convex/concave (see below) line bundles on $P$. Define the generating 
function of genus zero one-point gravitational Gromov--Witten invariants
on $X$ to be
\begin{equation} \label{e:1}
\begin{split}
 J_X (q,\hbar)
 := &e^{\frac{1}{\hbar} (t_0+ \sum_{i} p_i t_i)} \sum_{\beta} q^{\beta} 
    J_X(\beta)  \\
 := &e^{\frac{1}{\hbar} (t_0+ \sum_{i} p_i t_i)}  \sum_{\beta} q^{\beta}
      \ev^{\vir}_* \frac{1} {\hbar(\hbar-\psi)} \\
 = &e^{\frac{1}{\hbar} (t_0+ \sum_{i} p_i t_i)} \sum_{mk} T_m g^{mk}
  \sum_{\beta} q^{\beta} \int_{[\mbar_{0,1}(X,\beta)]^{\vir}} 
  \ev^* (T_k) \frac{1} {\hbar(\hbar-\psi)} 
\end{split}
\end{equation}
and for $E \to X$:
\begin{equation} \label{e:2}
\begin{split}
  J_X^E(q,\hbar) 
 := &e^{\frac{1}{\hbar} (t_0+ \sum_{i} p_i t_i)} \sum_{\beta} q^{\beta}
  	J_X^E(\beta)
 := e^{\frac{1}{\hbar} (t_0+ \sum_i p_i t_i)} 
    \ev^{\vir}_* \frac{\e(E_{\beta})} {\hbar(\hbar-\psi)},
\end{split}
\end{equation}
which will be called \emph{$J$-function} of $X$ and $E \to X$ respectively.
It is easy to see that
\[
  \int_X J_X^E =   \sum_{\beta} q^{\beta} \int_{[\mbar_{0,1}(X,\beta)]^{\vir}} 
  \frac{\e(E_{\beta})} {\hbar(\hbar-\psi)} 
  \ev^* e^{\frac{1}{\hbar} (t_0+ \sum_{i} p_i t_i)}.
\]

For a line bundle $L:= i^*(\so_{P}(l))$ on $X$ define 
\[
H^L_{\beta} := \prod_{k=0}^{\la c_1(L), \beta\ra} (c_1(L) +k \hbar) 
\]
for $L$ convex and
\[
H^L_{\beta} :=  \prod_{k=\la c_1(L), \beta\ra +1}^{-1} (c_1(L) +k \hbar) 
\]
for $L$ concave. Here a vector bundle $E \to X$ is called \emph{convex} if
for any stable map $f:C \to X$, $H^1(C, \pi_* f^*(E))=0$. $E$ is \emph{concave}
if $H^0(C, \pi_* f^*(E))=0$ for any $f$. Introduce another generating 0
of Gromov--Witten invariants on $X$ (modified by $H^L$'s)
\begin{equation} \label{e:iex}
  I^E_X(q,\hbar) := e^{\frac{1}{\hbar} (t_0+ \sum_{i} p_i t_i)} 
  \sum_{\beta} q^{\beta} J_X(\beta) \prod_j H^{L_j}_{\beta} 
\end{equation}

\begin{theorem} \label{t:1}
Let $X$ be a smooth projective variety embedded in $P=\prod \pp^{r_i}$
(see \emph{Notations.}) and $E=\oplus_j L_j \to X$ be the sum of line bundles 
$L_j$ which are the pull-backs of convex and concave line bundles on $P$ 
such that
\[
 c_1(T_X)-\sum_{L_j \text{convex}} c_1(L_j) +\sum_{L_j concave} c_1(L_{j})
\]
are non-negative. 
Then 
\[
 \int_X J_X^E(q,\hbar) \simeq \int_X I^E_X(q,\hbar),
\]
where $\simeq$ means equivalence up to a \emph{mirror transformation},
which is a special kind of change of variables described in \S~\ref{ss:4.4}.
\end{theorem}

In fact, there are many situations when the above mirror transformations
are unnecessary. The following theorem contains the main examples.

\begin{theorem} \label{t:2}
$\int_X J_X^E = \int_X I_X^E$ if
\begin{enumerate}
\item $E$ is concave and $\operatorname{rank}(E) \ge 2$. \label{concave}
\item $E$ is convex and $c_1(T_X)-\sum_{L_j \text{convex}} c_1(L_j)$ 
is Fano of index $\ge 2$. \label{convex}
\item Direct sum of the previous two cases.
\end{enumerate}
\end{theorem}

In the case when $E:= (\oplus_j L_{j})$, all $L_j$ are \emph{convex}, 
Theorem~\ref{t:1} has the following interpretation.

\begin{corollary}
Let $i_Y: Y \hookrightarrow X$ be the smooth zero locus of a section of 
$E$ (i.e.~a complete intersection in $X$).

\begin{list}{}{\setlength{\leftmargin}{0in}}
\item \textbf{1.} 
 \[  \int_X (i_Y)_* J_Y = \int_X J_X^E.
 \]
In particular, Theorem~1 and Theorem~2 relate the Gromov--Witten invariants 
of $Y$ to Gromov--Witten invariants of $X$.
\item \textbf{2.}
Suppose that $H^*(X)$ is generated by divisor classes, 
and $\rk(E)=1$, i.e.~$Y$ is a hypersurface. Then
all $n$-point gravitational Gromov--Witten invariants of $i_Y^* H^* (X)$ 
can be reconstructed from
one-point gravitational Gromov--Witten invariants of $X$.
\end{list}
\end{corollary}

\begin{proof}
The part 1 is nothing but the statement that for any $\omega \in H^*(X)$
\begin{equation} \label{e:ckl}
  \la [\mbar_{0,1}(Y,\beta)]^{\vir}, i_Y^*(\omega) \ra
  = \la [\mbar_{0,1}(X,\beta)]^{\vir} \cap \e(E_{\beta}), \omega \ra,
\end{equation}
which can be found in e.g.~\cite{CKL} and references therein. 
See also \S~\ref{ss:1.2}. The second part
is a corollary of the first part plus a reconstruction theorem proved in
\cite{LP}, which states that one can reconstructs $n$-point descendants
provided that $H^* (X)$ is generated by divisor classes and one-point 
descendants are known.
\end{proof}

\begin{remark} (Local mirror conjecture)
When $E$ is concave, 
there is also an interpretation of the Gromov--Witten invariants of $E
\to X$: the invariants of the total space of vector bundle $E \to X$. 
For example $\so(-1)\oplus \so(-1) \to \pp^1$, is the
``tubular neighborhood'' of $\pp^1$ embedded in a Calabi--Yau threefold 
with normal bundle $\so_{\pp^1}(-1)\oplus \so_{\pp^1}(-1)$.
\end{remark}

\begin{remark} \label{r:2}
Whenever $X$ carries a group action by $G$, one could carry out the
whole work to $H^*_G(X)$ (instead of $H^*(X)$) without any change. 
In this case, only localization of $\cc^*$-action on graph 
space is needed.
\end{remark}

Combining the quantum differential equation \cite{ABG1} and 
the above theorems, 
one can easily see the following interesting phenomenon. 
When $X$ is a toric variety this is a corollary of Givental's 
\emph{quantum Serre duality} theorem \cite{ABG2}. 

\begin{corollary} \label{c:qsd}
Let $E=\oplus_{j=1}^{\rk(E)} L_j$ be a concave bundle with $\rk(E) \ge 2$. 
Then\footnote{The restriction $\rk(E) \ge 2$ is shown unnecessary
in the joint work with A.~Bertram.}
\begin{equation} \label{e:qsd}
 \int_X \e(E) e^{\frac{1}{\hbar}(t_0 + p t)} J^E_X(q,\hbar)
 \sim \int_X \frac{1}{\e(E^{\vee})} e^{\frac{1}{\hbar}(t_0 + p t)} 
 J^{E^{\vee}}_X(q,\hbar),
\end{equation}
where $E^{\vee}$ is the dual vector bundle and $\sim$ means equivalence up  
to mirror transformations (and a factor of power series in $q$).
\end{corollary}

An interesting consequence of \eqref{e:qsd} is that one can prove the mirror
conjecture of convex bundles by concave bundles. For example the proof of
mirror conjecture in the case of quintic three-fold can be carried 
out\footnote{Here we are not very precise. \eqref{e:qsd} is valid only if 
$\rk(E) \ge 2$. However, we could use $E=\so(-5) \oplus \so(-1)$ on $\pp^5$
instead of $\so(-5)$ on $\pp^4$. The reason is that quintic three-fold can
be described as complete intersections in $\pp^5$ of the bundle $\so(1) \oplus
\so(5)$.}
by using $\so(-5)$ instead of $\so(5)$ on $\pp^4$. See \cite{ABG2} \S 5 for
an example. Here of course, this comes as a cyclic argument as we
have used the proof of mirror conjecture in convex case to deduce this
result. It is therefore desirable to have a direct proof of \eqref{e:qsd}
and refine the statement.
We plan to elaborate on this in a future paper (jointly with A.~Bertram).

\subsection{Relation to quantum cohomology} \label{ss:1.2}
Lefschetz hyperplane theorem (LHT) asserts that a projective smooth 
variety $X$ contains essential (co-)homological information of its hyperplane 
section $Y$. The quantum Lefschetz hyperplane theorem (QLHT) verifies this 
assertion in quantum cohomology. It was proposed by A.~Givental and 
formulated in the present form by B.~Kim \cite{BK}.

More precisely, our theorems state that the $J$-function of the complete 
intersection $Y$ (of classes $i_Y^* H^*(X)$) can be obtained from the
$J$-function of $X$ by multiplying suitable cohomology classes and
possibly a well-defined change of variables. 
The central roles of $J$-function in quantum cohomology theory is explained 
by the theory of \emph{quantum differential equation} developed by
Dijkgraaf and Givental \cite{ABG1}. It says, first of all, that
$J$-function is a flat section of the Dubrovin connection on the A-model
side, parallel to the Picard-Fuchs equation on the B-model side. This is 
related to the mirror symmetry discussed in the next subsection. 
Secondly, it gives a nice way to obtain the essential information of small 
quantum ring $QH^*(X)$ from $J$-function (\emph{quantum
$\mathcal{D}$-module}). For example, one can easily obtain the relations in 
quantum cohomology from $J$-function. Moreover, a result in \cite{LP} states 
that it is possible to reconstruct n-point descendants from 
one-point descendants when $H^*(X)$ is generated by divisor classes. 
Therefore our theorem even implicitly relates their big gravitational 
quantum cohomology algebras under this condition.

Many important special cases
of QLHT, including the celebrated quintic three-fold, have been worked out 
by A.~Givental (\cite{ABG2} and references therein), B.~Kim \cite{BK} and 
Liu--Lian--Yau \cite{LLY1,LLY2,LLY3}. In the case of quintic three-fold, 
the ambient space $X = \pp^4$ and $E=\so(5)$. Since it is very easy to compute
the Gromov--Witten invariants of $\pp^4$, QLHT is therefore the 
central part of the proof of mirror conjecture which will be discussed 
in the next subsection.

\subsection{Relation to mirror conjecture}
In a seminal paper \cite{CdGP} Candelas, de la Ossa, Green and Parkes applied
the mirror symmetry to the quintic three-fold and derived, in a 
string-theoretic way, the celebrated formula which predicts the number $n_d$ 
of rational curves on the quintic three-fold of any degrees.
This formula was then named \emph{mirror conjecture}, or mirror identity to 
distinguish itself from more fundamental physical principle of mirror symmetry.
This conjecture basically says that a generating function of 
$n_d$ is equivalent to a hypergeometric series up to a mirror transformation.
Their result soon stimulated a lot of mathematical work in enumerative
geometry. 
Among different groups working on the proof of their prediction,
there have been notably two different approaches.
One approach is trying to mathematically justify the string-theoretic mirror
symmetry and therefore obtain mirror conjecture as a corollary.
The other is to attack the enumerative consequence directly by developing
new mathematics inspired from physics. Gromov--Witten theory is partly
inspired by this second approach which we will give a brief 
discussion\footnote{The following is not meant to be a precise historical
account.}.

The first major progress came from M.~Kontsevich \cite{MK}. As is well known
in algebraic geometry, an enumerative problem can usually be formulated
as an intersection-theoretic one on suitable moduli spaces. Kontsevich
introduced the moduli space of stable maps and formulated the mirror conjecture
as follows. Let $\so(5)'_d:= \pi_* \ev^*(\so(5))$ be the vector bundle on 
$\mbar_{0,0}(\pp^4,d)$, then the enumerative problem was equivalent to computing
the integral 
\begin{equation} \label{e:mc1}
 N_d:=\int_{\mbar_{0,0}(\pp^4,d)} \e(\so(5)'_d),
\end{equation}
and $N_d$ can be related to $n_d$ by Aspinwall-Morrison Formula.
By using torus action on $\pp^4$ and fixed point localization method, he 
was able to reduce the integral \eqref{e:mc1} to summation of trees, 
but failed to complete the complicated combinatorial problem. 
Another (conceptual) drawback of this approach was that it did not
explain the presence of hypergeometric series.

Then came A.~Givental's proof followed by other approaches and generalizations
by Lian-Liu-Yau and Bertram. The new innovations include,
among other things, the introduction of equivariant quantum cohomology
and graph space. The quantum cohomology of Calabi-Yau manifold $X$ is
not semisimple, which makes the structure of quantum ring, like associativity
relation, not very useful in computing $n_d$. By introducing equivariant
quantum cohomology one produces a family of Frobenius structure over 
$H^*_G(pt)$ whose generic fibre carries semisimple Frobenius structure
while the special fibre $H^*(Y)$ does not.
Therefore $QH^*(Y)$ for Calabi-Yau manifold $Y$ may be considered as a 
limiting case of semisimple Frobenius manifold. This explains, in one way,
why the structure of quantum cohomology of Calabi-Yau manifolds play a role
in enumerative problem. 
On the other hand, to properly explain the presence of hypergeometric series, 
the graph space $\gbar_{0,0}(\pp^r,d)$ was introduced and was shown to have a 
natural birational morphism $u$ to the toric compactification space, or linear
sigma model, $\pp^r_d := \pp^{(r+1)d+r}$ (see \S~\ref{ss:gs} for details). 
It was earlier found by E.~Witten \cite{EW} and Givental, etc. (from different 
approaches) that some suitable correlators on $\pp^r_d$ actually produce the 
desired hypergeometric series.
However, neither $\gbar_{0,0}(\pp^r,d)$ nor $\pp^r_d$ is the right space to 
perform the integral \eqref{e:mc1}. The way to resolve this issue was to
identify $\mbar_{0,1}(\pp^r,d)$ as a fixed point component of $\cc^*$-action on 
$\gbar_{0,0}(\pp^r,d)$. One then uses the birational morphism
$u:\gbar_{0,0}(\pp^r,d) \to \pp^r_d$ to pass the above correlators from
$\pp^r_d$ to $\gbar_{0,0}(\pp^r,d)$ and then pass to $\mbar_{0,1}(\pp^r,d)$.
In the quintic three-fold case, $r=4$ and the correlator obtained from 
this procedure is
\begin{equation} \label{e:mc2}
 \begin{split}
  J_{\pp^4}^{\so(5)}(d) := 
   &\ev_* \frac{\e(\so(5)_d)}{\hbar(\hbar- \psi)} \\
  = & \hbar^{-2} \ev_* \e(\so(5)_d) + \hbar^{-3} \ev_* (\e(\so(5)_d) \psi)
	+\cdots. \\
 \end{split}
\end{equation}
The $\hbar^{-2}$ term in Laurent series expansion, when integrated over 
$\pp^4$, will be (see \cite{LLY1})
\[ \int_{\pp^4} \ev_* \e(\so(5)_d) = \int_{\mbar_{0,1}(\pp^4,d)} \e(\so(5)_d) 
  = \int_{\mbar_{0,0}(\pp^4,d)} \e(\so(5)'_d)
\]
which is exactly \eqref{e:mc1}.

There are now four approaches to mirror conjecture (known to us) by Givental,
Lian-Liu-Yau, Bertram, and Gathmann. The interested reader can find valuable 
information in \cite{ABG2,BK,LLY3,AB,AG} and references therein.

\begin{remark}
Compare \eqref{e:mc2} with \eqref{e:1}, one sees that our results can be
interpreted as a generalization of mirror conjecture.
\end{remark}


\begin{acknowledgements} 
My special thanks go to Aaron Bertram. The current proof is based on his work 
\cite{AB} and I benefit a lot from our collaboration. 
I am also thankful to Bumsig Kim\footnote{B.~Kim informed us that he had 
previously obtained a special case of Theorem~\ref{t:2} case~\ref{convex}.}, 
Rahul Pandharipande for numerous useful discussions.
\end{acknowledgements}

\begin{update} Gathmann has recently posted his proof of mirror theorem
\cite{AG2}, where he proved the mirror theorem in the case $E$ is a convex
line bundle (i.e.~when $Y$ is a very ample hypersurface). The relation
between his approach and ours is, roughly, the following.
While we tried to sweep the classes $e_{\nu}$ in Theorem~\ref{t:4}
``under the carpet'' by dimensional constraints, he explicitly studies
these relative classes and found a nice formula to relate these classes 
to ordinary Gromov--Witten classes.
%
\end{update}

\section{Graph space in Gromov--Witten theory}
\emph{The Picard number of $X$ is assumed to be one throughout the rest of 
the paper, to avoid complicated notations.
The generalization to the arbitrary Picard number is usually a 
matter of bookkeeping and is left to the reader.}

\subsection{Graph space and one-point invariants} \label{ss:gs}
The $n$-pointed \emph{graph space} of $X$ of degree $\beta$ is defined to 
be $\gbar_{0,n}(X,\beta):= \mbar_{0,n}(X \times \pp^1, (\beta,1))$,
where the degree $(\beta,1)$ is the element in $H_2(X)\oplus H_2(\pp^1)$. 
It is a compactification of the space of maps from parameterized 
$\pp^1$ to $X$. This space $\gbar_{0,n}(X,\beta)$ carries a natural 
$\cc^*$-action induced from the action on $\pp^1$.

When $X = \pp^r$ there is another (toric) compactification,  $\pp^r_d$ 
(\emph{linear sigma model}), of the space of parameterized map of degree $d$. 
It is constructed in the following way.
Consider the projective space of $(r+1)$-tuple of the degree $d$ (symmetric) 
binary forms of $(z_0:z_1)$. When there is no common factors of positive
degrees among these $(r+1)$-tuples, it represents a morphism from 
$\pp^1 \to \pp^r$. We may compactify it by simply allowing the common factors
and taking quotient by $\cc^*$-action. 
It is easy to see that this space is equal to $\pp^{(r+1)(d+1)-1}$. 
 
By construction $\gbar_{0,0}(\pp^r,d)$ is birationally isomorphic to
$\pp^r_d$. In fact, 

\begin{theorem} \label{t:3}
\emph{(Givental's Main Lemma \cite{ABG1})}
There exists a natural birational $\cc^*$-equivariant
morphism $u : \gbar_{0,0}(\pp^r,d) \to \pp^r_d$. 
\end{theorem}

This morphism $u$ can be described as follows. 
Consider a stable degree $(d,1)$ map $f : C\to \pp^1\times \pp^{n}$.
There exists a unique irreducible component $C_0\in C$ (called 
\emph{parameterized component}) such that $f |_{C_0}$
has degree $(d_0,1)$ where $d_0 \leq d$. The image $f (C_0)$ is the graph of
a map $\pp^1\to \pp^r$ of degree $d_0$. The map is given by the binary
forms $(p_0:\cdots:p_r)$ of degree $d_0$ with no common factors 
and determines the forms uniquely up to a non-zero constant factor. The curve 
$\overline{C \setminus C_0}$ has $s$ connected \emph{unparameterized 
components} which are mapped to $\pp^1\times \pp^r$ with
degrees $(d_1,0),...,(d_s,0),  d_1+...+d_s=d-d_0$, and the image of $i$-th
component is contained in the slice $(a_i :b_i)\times \pp^{n}$.
We put $u (C,f)=\prod_{i=1}^r (a_i z_0 -b_i z_1)^{d_i} 
(p_0:\cdots:p_r)$.    
For a detailed proof in algebro-geometric terms see \cite{LLY1} and 
\cite{BdPP}.

Of course we can also find simple birational models for $\gbar_{0,s}(\pp^r,d)$.
A particular useful one is $(\pp^1)^s \times \pp^r_{d}$. There is also
a morphism $u_s : \gbar_{0,s}(\pp^r,d) \to (\pp^1)^s \times \pp^r_{d}$.
The first factor is defined by the composition
\begin{equation*}
  \gbar_{0,s}(\pp^r,d) \overset{\prod \ev_i}{\longrightarrow} 
   (X\times \pp^1)^s \overset{\prod p_2}{\longrightarrow} (\pp^1)^s
\end{equation*}
and the second factor is the composition of the forgetful morphism
$\gbar_{0,s}(\pp^r,d) \to \gbar_{0,0}(\pp^r,d)$ and $u: \gbar_{0,0}(\pp^r,d)
\to \pp^r_d$ defined in Theorem~\ref{t:3}.
  
From the above description, it is clear that two spaces $\gbar_{0,0}(\pp^r,d)$
and $\pp^r_d$ differ on certain boundary strata. 
The \emph{comb type strata}, denoted $D_{\mu}$ where 
$\mu:=(d_0,d_1,\cdots,d_s), d_0+d_1+\cdots d_s=d$,
play an important role in our discussion. The
(domain) curve of a generic element in $D_{\mu}$ has one parameterized 
component and $s$ nodes, and $(d_0,1)\in H_2(\pp^r \times \pp^1)$ is the
degree of the parameterized component. These strata have substrata, called
\emph{hairy comb type strata}, which are obtained by further degenerating
the unparameterized components. 
Note that different permutations of $d_1,\cdots,d_s$ represent
the same $\mu$. For example $(d_0, d',d'')=(d_0,d'',d')$.
For such a stratum $D_{\mu}$ there is a finite birational morphism 
\begin{equation*}
  \tilde{D}_{\mu}:= \gbar_{0,s}(\pp^r,d_0) \times_{(\pp^r)^s}
  \prod_{m=1}^s \mbar_{0,1}(\pp^r,d_m) \to D_{\mu}.
\end{equation*}
To simplify our notation, we will denote $\tilde{D}_{\mu}$ also by $D_{\mu}$ 
henceforth. There are also useful morphisms from ${D}_{\mu}$ to simple spaces:
\begin{equation} \label{e:3}
  u_{\mu}: D_{\mu}= \gbar_{0,s}(\pp^r,d_0) \times_{(\pp^r)^s}
  \prod_{m=1}^s \mbar_{0,1}(\pp^r,d_m)	\to (\pp^1)^s \times \pp^r_{d_0},
\end{equation}
defined by the composition 
\begin{equation} \label{e:4} 
  {D_{\mu}} \overset{p_0}{\to} \gbar_{0,s}(\pp^r,d_0) 
	\overset{u_s}{\to} (\pp^1)^s \times \pp^r_{d_0}.
\end{equation}

Let $X \overset{i}{\hookrightarrow} \pp^r$ be an embedding described
earlier. For the notational convenience, we define the divisors $D_{\nu}
:= i_G^* D_{\mu}$ of $\gbar_{0,0}(X,\beta)$ (indexed by
$\nu:=(\beta_0,\beta_1, \cdots, \beta_s)$) for future reference. 
The relation is described in the following commutative diagram:
\begin{equation} \label{e:6}
 \xymatrix{\gbar_{0,0}(X,\beta) \ar[r]^{i_G} &\gbar_{0,0}(\pp^r,d) \ar[rr]^u
	&&\pp^r_d \\
 {D}_{\nu}  \ar[r]^{i_{\nu}} \ar[u]_{\varphi_{\nu}} 
	\ar@/_2pc/[rrr]_{u_\nu}
 &{D}_{\mu} \ar[rr]^{u_{\mu}} \ar[u]_{\varphi_{\mu}} 
   && (\pp^1)^s \times \prod \pp^r_{d_0} \ar[u]_{\psi_{\mu}}}
\end{equation}

The reason to introduce the graph space in Gromov--Witten theory is to
obtain one-point descendant invariants. The graph space
$\gbar_{0,0}(X,\beta)$ carries a $\cc^*$ action induced from the $\cc^*$
action on $\pp^1$, and $\mbar_{0,1}(X,\beta)$ is a fixed point component.
It is summarized in the following commutative diagram:

\begin{equation} \label{e:7}
\xymatrix{ \gbar_{0,0}(X,\beta) \ar[r]^{i_G} \ar@/^2pc/[rrr]^{u_X} 
    &\gbar_{0,0}(\pp^r,d) \ar[rr]^u &&{\pp^r_d} \\
  {\mbar_{0,1}(X,\beta)} \ar[u]^{j_X} \ar[r]^{i_M} \ar[d]_{\ev_X} 
	&\mbar_{0,1}(\pp^r,d) \ar[u]^{j_P} \ar[rr]^{\ev_P} &&P \ar[u]^t \\
  X \ar[urrr]_i }
\end{equation}
where the left upper square is a commutative diagram such that  
$\mbar_{0,1}(X,\beta)  \hookrightarrow \gbar_{0,0}(X,\beta)$ and 
$\mbar_{0,1}(P,d) \hookrightarrow \gbar_{0,0}(P,d)$ are fixed point components
of $\cc^*$-action on the graph spaces and $\mbar_{0,1}(X,\beta)$ is the only
fixed point component mapping to $\mbar_{0,1}(\pp^r,d)$ by $i_M$. The same
can be said about the right square: $\pp^r \hookrightarrow \pp^r_d$ is a
fixed point component such that $\mbar_{0,1}(\pp^r,d)$ is the only fixed
point component mapping to $\pp^r$. 

\begin{remark}
The above setting actually works for product of projective spaces 
$P:=\prod \pp^{r_i}$. There are only small changes in this adjustment.
First, all $d,\mu$ etc.~should stand for multi-index. Second,
the birational morphism in Theorem~\ref{t:3} should be replaced by
\begin{equation*} 
 u:\gbar_{0,0}(P,d) \to  \prod_i \pp^{r_i}_{d^i}. 
\end{equation*}
\end{remark}

\subsection{Virtual localization on graph space}
First recall the Graber-Pandharipande virtual localization formula 
\cite{GP}:

\begin{theorem} \label{t:vl}
Let $X$ be an algebraic scheme with a $\cc^*$ action and $\cc^*$-equivariant
perfect obstruction theory. Then the virtual localization formula holds:
\begin{equation*} 
[X]^{\vir}= j_* \sum \frac{[X_j]^{\vir}}{e(N_{X_j|X}^{\vir})}
\end{equation*}
in $A_*^{\cc^*}(X)\otimes \qq[\lambda, \frac{1}{\lambda}]$, where
$\lambda$ is the generator of the $A_*^{\cc^*}(pt)$.
\end{theorem}

An immediate consequence of this theorem is the 
\emph{correspondence of residues} (\cite{LLY3} Lemma~2.1 and \cite{AB}):

\begin{corollary} \label{c:1}
Suppose that $f:X_1 \to X_2$ is a $\cc^*$-equivariant map of two algebraic
schemes and $j_1: F_1 \to X_1$ and $j_2: F_2 \to X_2$ are two fixed point
components of $X_1$ and $X_2$ respectively, 
such that $F_1$ is the only fixed point component mapping
into $F_2$ as in the following commutative diagram:
\[
\begin{CD}
 X_1 @>f>> X_2 \\
 @Aj_1AA  @Aj_2AA \\
 F_1 @>{f|_{F_1}}>> F_2 .
\end{CD}
\]  
Then 
\begin{equation} \label{e:vl2}
{f|_{F_1}}_* \left(\frac{j_1^*(\omega) \cap [F_1]^{\vir}}
   {e(N^{\vir}_{F_1 |X_1})} \right) 
= \frac{j_2^{!} f_* (\omega \cap [X_1]^{\vir})}{e(N^{\vir}_{F_2 |X_2})}
\end{equation}
for any $\omega \in H^*_{\cc^*}(X_1)$.
\end{corollary}

Apply this result to our case: Let $X_1=\gbar_{0,0}(X,\beta)$, $X_2=
\gbar_{0,0}(P,d)$, $F_1= \mbar_{0,1}(X,\beta)$, $F_2= \mbar_{0,1}(P,d)$ and
$\omega = \e(E^G_{\beta})$, as displayed in the upper left square of
\eqref{e:7}. Since $ j_X^* \e(E^G_{\beta}) = \e(E_\beta)$, one has
\begin{equation} \label{e:vl3}
\begin{split}
 {i_M^{\vir}}_* \biggl(\frac{\e(E_\beta)}{\hbar(\hbar-\psi)}\biggr)
 =&\pd {i_M}_* \biggl(\frac{\e(E_\beta)}{\hbar(\hbar-\psi)} \cap 
  [\mbar_{0,1}(X,\beta)]^{\vir} \biggr) \\
 =&\frac{\pd j_P^! {i_G}_* \bigl(\e(E^G_{\beta}) \cap 
  [\gbar_{0,0}(X,\beta)]^{\vir} \bigr)} {\hbar(\hbar-\psi)} \\
 = &\frac{j_P^* {i_G^{\vir}}_* \e(E^G_{\beta})}{\hbar(\hbar-\psi)},
\end{split}
\end{equation}
where $\hbar$ is the generator of $H^*_{\cc^*}(pt)$ and
$\psi$ is the first chern class of the tautological line bundle $\mathcal{L}_1$
on $\mbar_{0,1}(X,\beta)$ or $\mbar_{0,1}(P,d)$.  The one small difference
between \eqref{e:vl2} and \eqref{e:vl3} is that $\gbar_{0,0}(P,d)$ and
$\mbar_{0,1}(P,d)$ are orbifolds and Poincar\'e duality makes sense there.

\begin{remark}
The functorial properties of virtual fundamental classes used in this article
can be found in \cite{LT} \cite{BF} \cite{KB}.
\end{remark}

\section{Decomposition of the virtual fundamental classes}

Recall that there is a birational morphism $u_{X,1}:=u_1\circ i_G$ 
(see \eqref{e:3}) from the universal curve $\gbar_{0,1}(X,\beta)$ of the 
graph space $\gbar_{0,0}(X,\beta)$ to $\pp^1 \times \pp^r_d$. 
Given a line bundle $L:=i^* \so_{\pp^r}(l)$ on $X$ one could produce two line 
bundles $\ev^* L$ and $u_{X,1}^* \so_{\pp^1 \times \pp^r_d}(dl,l)$ on 
$\gbar_{0,1}(X,\beta)$. 
These two line bundles are isomorphic on the open subset 
$U:=\gbar_{0,1}(X,\beta) \setminus \cup_{\nu=(\beta_0,\beta_1)} C_{\nu}$, 
where $C_{\nu}$ is the universal curve of the \emph{unparameterized} component 
over $D_{\nu}$.
The reason is that $C_{\nu}$ are exactly the exceptional divisors of $u_1$. 
More explicitly, there is a rational map $b: \pp^1\times \pp^r_d \to \pp^r$ 
such that $b$ is well-defined on $u(U)$ in the following commutative 
diagram:
\begin{equation*}
\xymatrix{\gbar_{0,1}(X,\beta) \ar[r] \ar[d]_{\ev_X} 
    &\gbar_{0,1}(\pp^r,d) \ar[r] \ar[d] &\pp^1\times \pp^r_d
  \ar@{-->}[ld]_{b}\\
    X \ar[r] & \pp^r}
\end{equation*}
It is easy to see that $b$ is a morphism on $u(U)$, of degree $d$ in the first 
factor and linear in the second. Namely 
\[
  b|_{u(U)}^* \so_{\pp^r} (l) =\so_{\pp^1 \times \pp^r_d}(dl,l). 
\]
Therefore 
\[
 \ev_X^*(L)  =u_{X,1}^*(\so_{\pp^1\times \pp^r_d}(dl,l))\otimes 
  \so \Bigl( \sum_{\nu=(\beta_0,\beta_1)} c_{\nu} C_{\nu} \Bigr) .
\]
It is easy to see that $c_{\nu}$ is negative as $\ev_X^*(L) \hookrightarrow
\so_{\pp^1\times \pp^r_d}(dl,l))$ by restricting a section of $\ev_X^*(L)$ to
$U$ and then extending this section on $\pp^r_d$ by Hartog's lemma. The actual 
coefficients\footnote{In fact, we will only need to know $c_{\nu}$ is 
negative.} $c_{\nu}$ can be determined by the following observation.
The degree of $f^*(L)$ on the \emph{unparameterized} component of the
universal curve over the divisor $D_{\nu}$ (at a generic point) 
has degree $\la c_1(L), \beta_1\ra $. This leads
to $c_{\nu}= - \la c_1(L), \beta_1\ra $. 
The same argument applies to $n$-pointed graph space:

\begin{lemma} \label{l:1}
On the universal curve $C$ over the graph space $\gbar_{0,n}(X,\beta)$
\begin{equation*} 
\begin{split}
 & \ev_{1,X}^*(L) \\
  =& u_{X,n}^*(\so_{\pp^1\times (\pp^1)^n \times \pp^r_d}
  (dl,0,\cdots,0,l)) \otimes \so \Bigl( \sum_{\nu=(\beta_0,\beta_1)} 
  -\la c_1(L), \beta_1\ra C_{\nu} \Bigr),
\end{split}
\end{equation*}
where $C_{\nu}$ is the universal curve of the unparameterized component over 
$D_{\nu}$. 
\end{lemma}

\begin{corollary} \label{c:2}
For $L$ convex let 
\begin{gather*} 
  L_{\beta}^G:=R^0 \pi_{n+1 *} \ev_{n+1}^*(L_l)\\
  F_{\beta}^G := u_X^* \Bigl(H^0(\pp^1,\so(dl)) \otimes
   \so_{(\pp^1)^n \times \pp^r_d}(0,\cdots,0,l)\Bigr). 
\end{gather*}
One has equivariant maps of vector bundles on graph space 
$\gbar_{0,n}(X,{\beta})$
\begin{equation*} 
  \sigma_0 : L_{\beta}^G \to  F_{\beta}^G.
\end{equation*}
For $L$ concave let 
\begin{subequations} 
\begin{gather}
  L_{\beta}^G := R^1 \pi_* e^*(L_l (-x_1-\cdots-x_n)) \label{e:11a}\\
  F_{\beta}^G := u_X^* \Bigl(H^1(\pp^1,\so(d l)(-\chi_1 -\cdots-\chi_n)) 
   \otimes \so_{(\pp^1)^n \times \pp^r_d}(0,l) \Bigr).  \label{e:11b}
\end{gather}
\end{subequations} 
One has
\begin{equation*} 
  \sigma_1 : F_{\beta}^G \to L_{\beta}^G.
\end{equation*}
Here in \eqref{e:11b} $\chi_m: (\pp^1)^n \otimes \pp^r_d \to \pp^1$ are
the projections to the $m$-th factor of $\pp^1$ (considered as  
``marked points''). In equation \eqref{e:11a}, $x_m$ are the marked points
on the universal curve. 
\end{corollary}

\begin{proof} 
It is clear that the vector bundles in this lemma are the push-forwards of
two line bundles on the universal curve $C$ considered in Lemma~\ref{l:1}.
When the coefficient $c_{\nu}$ is negative, one has the following inclusion
\begin{equation} \label{e:13}
 \ev_X^*(L) \hookrightarrow 
 u_{X,n}^*\so_{\pp^1\times (\pp^1)^n \times \pp^r_d}(dl, 0, d).
\end{equation}
$\sigma_0$ is then obtained by pushing-forward the above inclusion of line 
bundles to $\gbar_{0,n}(X,\beta)$. 

$\sigma_1$ can be obtained in a similar way. When $L$ is concave, 
$-\la c_1(L), \beta_1 \ra$ is positive so that the arrow of \eqref{e:13} 
is reversed. 
\end{proof}

The equation \eqref{e:ckl} indicates that the QLHT boils down to the study
of the top chern class $\e(E_{\beta})$ on $\mbar_{0,1}(X,\beta)$, which is
in turn the pull-back $j_X^* (\e(E^G_{\beta}))$ of top chern class from
$\gbar_{0,0}(X,\beta)$. From the above discussion
\[
  \e(E^G_{\beta}) = u_X^* \e(F_{\beta}^G) + \text{boundary terms},
\]
where the boundary terms are supported on the comb type strata $D_{\nu}$.
Notice that the ``main term'' $\e(F_{\beta}^G)$ is the (equivariant) top
chern class of direct sum of line bundles on the projective space and can
therefore be easily computed. In fact, we will see
in the next section that this part gives rise to the factors $H^E_{\beta}$ in
\eqref{e:iex}. It remains to have a close look of the boundary terms. A key
observation of \cite{AB} is that boundary contributions can be explicitly
computed using MacPherson's graph construction for vector bundle morphisms.

\begin{proposition} \label{p:0}
When $X=\pp^r$ and $E=\so_{\pp^r}(l)$ a line bundle, 
\begin{equation*} 
 \e (\so_{\pp^r}(l)^G_{d})=
  u^* \prod_{k=0}^{dl} (lH+k \hbar)
  + \sum_{\mu} \frac{1}{s!} {\varphi_{\mu}}_* \bigl(e_{\mu} \cup 
 u_{\mu}^* \prod_{k=0}^{d_0 l} (l H +k \hbar)\bigr),
\end{equation*}
for $l$ \emph{positive} and
\begin{equation*} 
  \e (\so_{\pp^r}(l)^G_{d})
  =u^* \prod_{k=d l+1}^{-1} (l H+k \hbar) 
  + \sum_{\mu} \frac{1}{s!} {\varphi_{\mu}}_* \Bigl(e_{\mu} \cup 
   u_{\mu}^* \xi_1 \cdots \xi_s \prod_{k=d_0 l+1}^{-1} (l H +k \hbar) \Bigr),
\end{equation*}
for $l$ negative. Here $s$ is the number of unparameterized components for
a generic curve over $D_{\nu}$, $H$ is the hyperplane class of 
$\pp^{r}_{d}$ and $\xi_m$ is the point class of $m$-th $\pp^1$ in 
$(\pp^1)^s \times \pp^r_d$. The equivariant class $e_{\mu}$ is defined
in \eqref{e:30}. 
\end{proposition}

\begin{theorem} \label{t:4}
Let $E=\oplus_j L_j$, $L_j=i^*(\so_{\pp^r}(l_j))$.
The virtual fundamental class $\e(E^G_{\beta}) \cap 
[\gbar_{0,0}(X,\beta)]^{\vir}$ 
decomposes as follows: 

For \emph{convex} bundle $E$
\begin{equation} \label{e:37}
\begin{split}
  &\e  (E^G_{\beta}) \cap [\gbar_{0,0}(X,\beta)]^{\vir} \\
  = & u_X^* \prod_j \prod_{k_j=0}^{\la c_1(L_j), \beta \ra} (l_j H+k_j \hbar) 
  \cap [\gbar_{0,0}(X,\beta)]^{\vir} \\
  + &\sum_{\nu} \frac{1}{s!} 
  	{\varphi_{\nu}}_* \Bigl( \bigl(e_{\nu} \cup u_{\nu}^* \prod_j 
  \prod_{k_j=0}^{\la c_1(L_j), \beta_0\ra} (l_j H +k_j \hbar) \bigr) \cap
	[{D}_{\nu}]^{\vir} \Bigr). 
\end{split}
\end{equation}

For \emph{concave} $E$
\begin{equation*} 
\begin{split}
  &\e  (E^G_{\beta}) \cap [\gbar_{0,0}(X,\beta)]^{\vir}  \\
  = & u_X^* \prod_j \prod_{k_j=\la c_1(L_j),\beta \ra +1}^{-1} 
    (l_j H+k_j \hbar) \cap [\gbar_{0,0}(X,\beta)]^{\vir}\\
  + &\sum_{\nu} \frac{1}{s!} 
  {\varphi_{\nu}}_* \Bigl( \bigl(e_{\nu} \cup u_{\nu}^* \prod_j
 \prod_{k_j=\la c_1(L_j), \beta_0\ra +1}^{-1} \xi_1 \cdots \xi_s
   (l_j H +k_j \hbar) \bigr) \cap [{D}_{\nu}]^{\vir} \Bigr).
\end{split}
\end{equation*}
In the case $\rk(E)=1$, $e_{\nu}:=i_G^* e_{\mu}$. In general, it is defined
inductively.
\end{theorem}

\begin{proof}
It is easy to see that the theorem follows from the above proposition
by the formula $\e (E^G_{\beta})=\prod_j \e ((L_j)^G_{\beta})$ and 
the excess intersection theory. Note that $\e ((L_j)^G_{\beta})=
i_G^* \e(\so_{\pp^r}(l_j)^G_d)$.
\end{proof}

\begin{proof}  (of Proposition~\ref{p:0})

\textbf{I.} (\emph{Convex case}) 
The convex case of the proposition is Lemma~4.4 of \cite{AB}. For the  
convenience of the reader and future references, we reproduce Bertram's proof.

As remarked above, it is clear from Lemma~\ref{l:1} that the equivariant 
virtual class has its ``main'' contribution from 
\[
  \e (u^* \so_{\pp^1\times \pp^r_d}(ld,l)) 
   = \Bigl(u^* \prod_{k=0}^{ld} (lH + k\hbar)\Bigr) .
\]
It is also clear that other contributions to $\e (\so_{\pp^r}(l)^G_{d})$ come
from the boundary strata $\mu$ (and possibly its substrata). It remains to
study the boundary terms. This can be done by MacPherson's graph construction
of vector bundle morphisms.

Let $\E_d:=\so_{\pp^r}(l)^G_{d}$ and 
$\F_d := H^0(\pp^1, \so(dl))\otimes u^* \so_{\pp^r_d} (l)$ be vector bundles on
$\gbar_{0,0}(\pp^r,d)$. By Corollary~\ref{c:2} there is a homomorphism 
$\sigma_0: \E_{d} \to \F_{d}$ of vector bundles of the same rank $dl+1$.
Let $G:=\text{Grass}_{dl+1}(\E_{d}\oplus \F_{d})$ be the Grassmann bundle over 
$\gbar_{0,0}(\pp^r,d)$, with universal bundle $\zeta \to G$ of rank $dl+1$.
There is a canonical embedding
\begin{equation*} 
   \Phi: \gbar_{0,0}(\pp^r,d) \times  \mathbb{A}^1 \to G \times \pp^1
\end{equation*}
taking $(z, \lambda)$ to $(\text{graph of } \lambda \sigma_0(z), (1:\lambda))$.
Let $W$ be the closure of the image of $\Phi$ and let
\begin{equation*} 
   W_{\infty} = i_{\infty}^*[W] = \sum m_{\delta} [V_{\delta}] 
\end{equation*}
be a cycle in $G$ of dimension equal to $\dim(\gbar_{0,0}(\pp^r,d))$, where
$m_{\delta}$ is the multiplicity of $V_{\delta}$. Let 
\[
 \eta_{\delta} : V_{\delta} \subset G \to  Z_{\delta}
  \subset  \gbar_{0,0}(\pp^r,d)
\]
be the map induced by projection, and $Z_{\delta}$ be the image of 
$\eta_{\delta}$. 
The philosophy of the graph construction is that different components
$Z_{\delta}$ are responsible for different types of degeneration of $\sigma_0$.
Since $\sigma_0$ is generically of full rank, there is one component 
$V_0 \sim Z_0 \cong \gbar_{0,0}(\pp^r,d)$ in $W_{\infty}$ with multiplicity
one such that $Z_0$ embeds in $G$ via the fibre of $\F_{d}$.
Other $V_{\delta}$ map to proper sub-varieties $Z_{\delta}$ of 
$\gbar_{0,0}(\pp^r,d)$. It follows that (\cite{WF}, Example 18.1.6)
\begin{equation} \label{e:19}
 \e(\E_{d})
 = u^* \prod_{k=0}^{dl+1}(lH+k\hbar) 
  + \sum_{\delta} m_{\delta} {\eta_{\delta}}_* \e(\zeta). 
\end{equation}

To calculate the contributions from $V_{\delta}$, further study on the behavior
of $\sigma_0$ on boundary strata $D_{\mu}$ is needed. Let $f:C \to \pp^r$
be a generic stable map in strata $\mu=(d_0,\cdots,d_r)$.
$\sigma_0|_{D_{\mu}}$ can be described fibrewisely (at a generic point) as:
\begin{equation} \label{e:20}
  H^0 (C,{f}^* \so_{\pp^r}(l)) \to H^0 (C_0, f^*\so_{\pp^r}(l) |_{C_0}) 
  \cong H^0(\pp^1,\so(d_0 l)) \to H^0 (\pp^1, \so(dl)),
\end{equation}
where $C_0\cong \pp^1$ is the parameterized component of $C$.
The first map is simply the restriction and the last map is defined by 
multiplying a factor 
$\prod_{m=1}^s (a_m z_0- b_m z_1)^{l d_m}$, where
$(a_m:b_m)$ are the nodal points on the parameterized $\pp^1$ and $(z_0:z_1)$
the homogeneous coordinates on parameterized $\pp^1$. It follows that 
$\sigma_0$ drops ranks only on the comb strata of types $\mu$ described in the
previous section. Any hairy comb substrata obtained from $\mu$ by further 
degenerating the unparameterized components will not further reduce the rank 
of $\sigma_0$. Moreover, it has the following \emph{transversality} property:
$\sigma_0$ has generic corank $n_1$ and $n_2$ along $D_{\mu_1}$ and 
$D_{\mu_2}$ respectively. Then
$\sigma_0$ has generic corank $n_1+n_2$ along the intersection of two strata.

The above fibrewise description actually holds globally. 
Namely, on the strata $D_{\mu}$, $\sigma_0$ is the following composition 
\begin{equation} \label{e:21}
  \varphi_{\mu}^* \sigma_0  :\varphi_{\mu}^* \E_d
  \to  p_0^* \E_{d_0} \overset{p_0^* \sigma_0} {\to} 
  p_0^* \F_{d_0} \to \varphi_{\mu}^* \F_{d},
\end{equation}
where $p_0 : {D}_{\mu} \to \gbar_{0,s}(\pp^r,d_0)$ is the projection 
(see \eqref{e:4}), $\E_{d_0}$ and $\F_{d_0}$ on $\gbar_{0,s}(\pp^r,d_0)$ 
are defined in Corollary~\ref{c:2}.
The last map in \eqref{e:21} is the push-forward (along the first $\pp^1$) 
of the following map 
\begin{equation*}
  \so_{\pp^1\times (\pp^1)^s \times \pp^r_{d_0}}
	(d_0 l, 0,\cdots,0,l) 
  \hookrightarrow \so_{\pp^1\times (\pp^1)^s \times \pp^r_{d_0}}
        (dl,d_1 l, \cdots,d_s l,l).
\end{equation*}

One can now apply the above study of degeneration type of $\sigma_0$
to the graph construction. Since the rank of 
$\sigma_0$ decreases only on the strata $D_{\mu}$, \eqref{e:19} becomes
\begin{equation} \label{e:23}
 \e(\E_{d}) = u^* \prod_{k=0}^{dl}(lH+k\hbar) +\sum_{\mu} m_{\mu} S_{\mu}.
\end{equation}
Namely, the only $Z_{\delta}$ in \eqref{e:19} are $D_{\mu}$ and 
$S_{\mu}= {\eta_{\mu}}_* \e(\zeta)$. To write down $S_{\mu}$ explicitly
in terms of characteristic classes in $\sigma_0 |_{D_{\mu}}: 
\varphi_{\mu}^* \E_{d} \to \varphi_{\mu}^* \F_{d}$ 
one would need a filtration of the above setting.
Let us first deal with the simplest case when $\mu=(d_0,d_1)$, 
i.e.~$D_{\mu}$ is a divisor. The universal curve on 
${D}_{\mu}=\gbar_{0,1}(\pp^r,d_0) \times_{\pp^r} \mbar_{0,1}(\pp^r,d_1)$ 
(generically) consists of one parameterized and one unparameterized components.
The kernel of $\sigma_0$ can be identified with $p_1^* \E_{d_1}^1$, where 
$p_1: D_{\mu} \to \mbar_{0,1}(\pp^r,d_1)$ and $\E^1_{d_1}$ is the kernel
of the evaluation morphism $e$ of bundles on $\mbar_{0,1}(\pp^r,d_1)$
\begin{equation} \label{e:24}
  0 \to \E_{d_1}^1  \to \pi_{2 *} \ev_2^* \so_{\pp^r}(l) \overset{e}{\to} 
   \ev_1^*  \so_{\pp^r}(l) \to  0,
\end{equation}
where $\ev_2: \mbar_{0,2}(\pp^r,d_1) \to \pp^r$ and
$\pi_2:\mbar_{0,2}(\pp^r,d_1) \to \mbar_{0,1}(\pp^r,d_1)$ forgets the
second marked point. The kernel $\E^1_{d_1}$ can be further filtered by the
order of zeros of $e$:
\begin{equation} \label{e:25}
 0=p_1^* \E^{d_1 l+1}_{d_1} \subset p_1^* \E^{d_1 l}_{d_1} \subset \cdots 
 \subset p_1^* \E^1_{d_1} \subset \varphi_{\mu}^* \E_{d}
\end{equation}
where $\E^k_{d_1}$ consists of those sections of $\ev_2^* \so_{\pp^r}(l)$
which vanishes at least to the $k$-th order at the marking. Similarly we
can filter $\varphi_{\mu}^* \F_{d}$ by the span of the image of 
$\sigma_0|_{kD_{\mu}}$
on $\pp^1 \times \pp^r_{d_0}$:
\begin{equation*} 
 \F^k_{d_0} := H^0(\pp^1,\so(d_0 l+k-1))
  \otimes \so_{\pp^1 \times \pp^r_{d_0}}(k-1,l)
\end{equation*} 
such that
\[
  p_1^* \F_{d_0}= u_{\mu}^* \F_{d_0}^1 \subset \cdots \subset 
  u_{\mu}^* \F^{d_1 l}_{\mu} \subset \varphi_{\mu}^* \F_{d}
\]
is a filtration of $\varphi_{\mu}^* \F_{d}$ on $D_{\mu}$.
Now for each $D_{\mu}$ there are $d_1 l$ components 
$V^k_{\mu}\subset W_{\infty}$ ($\eta_{\mu}: V_{\mu} \to D_{\mu}$) because
the infinitesimal property of $\sigma_0$ along $D_{\mu}$. Each of $V^k_{\mu}$
is a birational image of a $\pp^1$-bundle over $D_{\mu}$. By a local
computation (\cite{AB}), $V^k_{\mu}$ has multiplicity\footnote{The explicit 
multiplicity is actually irrelevant to our result.} $k$, and the tautological
bundle on $V^k_{\mu}$ can be expressed in terms of the filtration of 
$\E_{d}$ and $\F_{d}$:
\begin{equation*} 
 0 \to p_0^* u_{\mu}^*\F^k_{d_0} \oplus p_1^* \E^{k+1}_{d_1} \to 
  \zeta_{V^k_{\mu}} \to \so(-1) \to 0.
\end{equation*}
The contribution $S_{\mu}$ from the boundary strata $\mu=(d_0,d_1)$ can
be therefore written as
\begin{equation*} 
  S_{\mu} = {\varphi_{\mu}}_* \bigl(e_{\mu} \cup 
  \e(u_{\mu}^*\F^1_{d_0}) \bigr)  
 ={\varphi_{\mu}}_* \bigl( e_{\mu} \cup u_{\mu}^* \prod_{k=0}^{d_0 l} 
  (lH + k \hbar) \bigr) 
\end{equation*}
where 
\begin{equation*} 
 e_{\mu}=\sum_{k=1}^{d_1 l} (-k) p_{d_1}^* 
  \e(\E_{d_0}^{k+1}) \cup u_{\mu}^* \e(\F^k_{d_0}/\F^1_{d_0}).
\end{equation*} 
This is exactly what we are looking for.

When $\mu$ is not a divisorial stratum we may use the above transversality
property (of $\sigma_0$ concerning the intersection of two boundary
strata). Since every $D_{\mu}$ is the intersection of divisors,
and the corank of $\sigma_0$ at intersection is equal 
to the sum of the coranks generically, we can then obtain $S_{\mu}$ for 
general $\mu$ inductively. Let $\mu_1=(d_0,d_1)$ and $V_{{\mu}_1}^k 
\subset W_{\infty}$ be the $\pp^1$-bundle over $D_{\mu_1}$. Now
apply graph construction to the following vector bundle morphism on 
(pulled back to) $V_{{\mu}_1}^k$
\begin{equation*} 
  \varphi_{\mu}^* \E_{d}/ p_1^* \E_{d_1}^1 \cong p_0^* \E_{d_0}
  \overset{p_0^* \sigma_0}{\to} p_0^* \F_{d_0} = u_{\mu}^* \F^1_{d_0}.
\end{equation*}
(Some obvious pull-backs will be omitted.)
Then the components of $W'_{\infty}$ over $V_{\mu_1}^k$ obtained from
this construction map birationally to the components of $W_{\infty}$ over
$\gbar_{0,0}(\pp^r,d)$. For example,
components of $V^m_{\mu'} \subset W'_{\infty}$ corresponding to boundary 
strata $\mu'=(d_0 -d',d')$ of $\gbar_{0,1}(\pp^r,d_0)$ 
map to the components of $V^m_{\mu} \subset W_{\infty}$
over the boundary strata $\mu=(d_0 -d', d',d_1)$ on
$\gbar_{0,0}(\pp^r,d)$. This implies that $V_{\mu}$ for any $\mu=(d_0,
d_1,\cdots, d_s)$ with a fixed ordering of $d_1,\cdots,d_s$,
$V^k_{\mu}$ are (birationally) towers of $\pp^1$-bundles over $D_{\mu}$
(by first doing $\mu_1=(d_0+\cdots+d_{s-1},d_s)$ then 
$\mu_2=(d_0+\cdots+d_{s-2},d_{s-1},d_s)$, etc.).  
An explicit expression of $e_{\mu}$ can therefore be obtained:
\begin{equation} \label{e:30} 
e_{\mu} = \prod_{m=1}^s \sum_{k_m=1}^{d_m l} (-k_m) p_{d_m}^* 
  \e(\E^{k_m+1}_{d_m}) \cup u_{\mu}^* \e(\F_{\Sigma_{m}}^{k_m}/
	\F^1_{\Sigma_{m}})
\end{equation}
where $\Sigma_{m}:= \sum_{a=0}^{m-1} d_a$ and $\E^k_{d_m}$ is defined 
to be the filtration of kernel on the $m$-th unparameterized
component, similar to that defined in \eqref{e:25} and $\F_{\Sigma_m}^1 
\hookrightarrow \F^{k_m}_{\Sigma_m}$ is the push-forward of the inclusion of
line bundles on $\pp^1 \times (\pp^1)^s \times P_{d_0}$ to 
$(\pp^1)^s \times P_{d_0}$:
\begin{equation*}
\begin{split}
  \so&(\sum_{a=1}^{m-1} d_a l, d_1 l, \cdots,d_{m-1} l, 0,\cdots,0,l) \\ 
  \to \so&(\sum_{a=1}^{m-1} d_a l +k_m-1, d_1 l, \cdots, d_{i-1} l,k_m-1, 0,
   \cdots,0,l).
\end{split}
\end{equation*}
This completes our proof of the convex case.

\textbf{II.} (\emph{Concave case}) 
The proof of the concave case can in general be carried out in a similar
way. However, some crucial modifications will be necessary. 
Now $\E_{d}=R^1 \pi_* \ev^* \so_{\pp^r}(l)$ and $\F_{d}=
H^1(\pp^1,\so(dl)) \otimes u^* \so_{\pp^r_d}(l)$ (with $l$ negative).
Corollary~\ref{c:1} guarantees that there is a bundle map $\sigma_1$ between
$\F_{d}$ and $\E_{d}$. We may apply Serre duality and find
$\sigma_1^* : (\E_{d})^* \to (\F_{d})^*$. 
Carry out the graph construction to $\sigma_1^*$. The equation \eqref{e:20}
(on $D_{(d_0, \cdots, d_s)}$) should be replaced by 
\begin{multline*} 
 H^0(C,\omega_C \otimes f^* L^{-1}) \to 
  H^0(C_0,\omega_{C_0}(\chi_1+\cdots+\chi_s) \otimes f^* L^{-1}|_{C_0}) \\
   \cong H^0(\pp^1, \omega_{\pp^1} (-d_0 l)\otimes \so(\chi_1+\cdots+\chi_s))
  	\to H^0(\pp^1,\omega_{\pp^1} (-dl)).
\end{multline*}
Here $\chi_1,\cdots,\chi_s$ are the nodal points on the parameterized $\pp^1$ 
(see Corollary~\ref{c:1}) and $\omega_C$ is the dualizing sheaf of $C$.

A similar modification to \eqref{e:23} should also take place:
\[
  \e(E_{d}) = u_X^* \prod_{k= \la c_1(L),d \ra +1}^{-1} (lH+k \hbar) 
   +\sum_{\mu} m_{\mu} S_{\mu}.
\]

For the filtration of the kernel of $\sigma_1^*$, we may use the following 
exact sequence (see \eqref{e:24})
\begin{equation*} \label{e:34}
 0 \to \E^{1}_{d_1} \to \pi_* (\ev_2^* \so_{\pp^r}(l) \otimes \omega(x_1)) 
 \overset{\operatorname{res}}{\to} \ev_1^* \so_{\pp^r}(l) \to 0,
\end{equation*}
where $\operatorname{res}$ is the residue at $x_1$. Again we can further
filter $\E^{1}_{d_1}$ by the order of zeros of $\operatorname{res}$,
(as did in \eqref{e:25}). Similar filtration can be defined on 
$\varphi_{\mu}^* \F_{d}$. Namely
\begin{equation*} \label{e:35}
\begin{split}
 \F^{k_m}_{d_0}=  &\omega_{\pp^1} \bigl( x_1+\cdots+x_s+
  \sum_{a=1}^{m-1} d_a (-l) +k_m-1 \bigr) \\
  &\otimes \so_{(\pp^1)^s \times P_{d_0}} \bigl(d_1 (-l), \cdots,d_{m-1} (-l),
  k_m-1, 0,\cdots,0,(-l)\bigr).
\end{split}
\end{equation*}

Now a similar computation leads to
\begin{equation} \label{e:36}
\begin{split}
  S_{\mu} = &{\varphi_{\mu}}_* \bigl( e_{\mu} \cup \e(u_{\mu}^*\F^1_{d_0}) 
   \bigr) \\
 =&{\varphi_{\mu}}_* \bigl(e_{\mu} \cup u_{\mu}^* \xi_1 \cdots \xi_s 
 \prod_{k=d_0 l +1}^{-1} (lH + k \hbar) \bigr).
\end{split}
\end{equation}
The rest is straightforward and is left to the reader.
\end{proof}

\section{Conclusion of the proofs} \label{s:4}

\subsection{Main contribution term}\footnote{Here 
instead of going through the diagrams \eqref{e:41} \eqref{e:43}, it might
be possible to proceed by another (equivalent) way using Givental's double 
construction formula \cite{ABG1}, which reads 
\[ \mathcal{G}:=
 \sum_{\beta} q^{\beta}\int_{[\gbar_{0,0}(X,\beta)]^{\vir}} e^{ P t} 
  \e (E_{\beta}) =\int_X J_X^E(q e^{\hbar t}, \hbar)
 e^{pt} J_X^E(q, -\hbar),
\]
where $P=u_X^*(\so_{\pp^{r}_{d}}(1))$.
This should give us $J_X^E(q,\hbar,t_0,t)$ (see \eqref{e:61}) 
from formulas in Theorem~\ref{t:4}.}
Set $\ev = i\circ \ev_X$ and 
\[
 J_X^E(\beta)=J_X^E(\beta, \text{main}) + \sum_{\nu} J_X^E(\beta,\nu).
\]
In order to show that $\int_X J_X^E \simeq \int_X I_X^E$ it is 
sufficient to show $i_* J_X^E \simeq i_* I_X^E$ because
\[
  \int_X e^{pt} \omega = \int_{\pp^r} e^{ht} i_* \omega
\]
by the projection formula.
 
Recall that our goal is to compute $\ev^{\vir}_* \dfrac{\e(E_{\beta})}
{\hbar (\hbar-\psi)}$. From \eqref{e:7} we have 
\begin{equation} \label{e:41}
\xymatrix{\gbar_{0,0}(X,\beta) \ar[rr]^{u_X} && \pp^r_d \\
 \mbar_{0,1}(X,\beta) \ar[u]_{j_X} \ar[rr]^{\ev}  &&\pp^r \ar[u]_{t} }
\end{equation}
We can easily get, by correspondence of residues \eqref{e:vl2}, the 
``main term'':
\begin{equation*} \label{e:42}
\begin{split}
 i_* J_X^E(\beta,\text{main})= 
   &\ev^{\vir}_* \left( \frac{j_X^* u_X^* \prod_j \prod_{k_j} 
  (l_j H+k_j \hbar)} {\hbar(\hbar-\psi)} \right) \\
 = &\frac{t^* {u^{\vir}_X}_* u_X^*\prod_j \prod_{k_j} (l_j H+k_j \hbar)}
 {e (N_{\pp^r|\pp^r_d})} \\
 = &\frac{t^* {u^{\vir}_X}_* 1} {e (N_{\pp^r|\pp^r_d})} \cup 
	t^*\prod_j \prod_{k_j} (l_j H+k_j \hbar)\\
 = &\ev^{\vir}_* (\frac{1}{\hbar (\hbar-\psi)}) \cup \prod_j \prod_{k_j}
   (l_j h +k_j \hbar) \\
 = &i_* I_X^E.
\end{split}
\end{equation*}
Thus $I_X^E$ is really the main term of $J_X^E$. We will
see that the boundary terms are of special forms and can be 
taken care of by a change of variables due to the \emph{non-negativity}
condition (on the tangent bundles of complete intersection $Y$ in $X$) 
stated in Theorem~\ref{t:1}.

\subsection{Boundary terms and dimension counting} \label{ss:boundary}
Before we start our discussion, we should remark that the term ``dimension''
here means \emph{virtual dimension}. 

Consider the commutative diagram (see \eqref{e:7})
\begin{equation} \label{e:43}
  \xymatrix{\gbar_{0,0}(X,\beta) \ar[r]^{u_X}  &\pp^r_d \\
   {D}_{\nu} \ar[r]^{u_{\nu}} \ar[u]^{\varphi_{\nu}}  
   &(\pp^1)^s \times \pp^r_{d_0} \ar[u]^{\psi_{\nu}} && \pp^r 
	\ar[ll]_{t_{\nu}}  \ar[llu]_{t} }_.
\end{equation}
In the convex case (see \eqref{e:37}),
\begin{equation} \label{e:44}
\begin{split}
  &i_*J_X^E(\beta,\nu) \\
 = & \ev^{\vir}_* \biggl( \frac{j_X^* \sum_{\nu} {\varphi_{\nu}}_* 
   \frac{1}{s!} \Bigl( e_{\nu} \cup u_{\nu}^* \prod_j \prod_{k_j=0}^{d_0 l_j}
  (l_j H +k_j \hbar) \cap [D_{\nu}]^{\vir} \Bigr)} 
  {\hbar (\hbar-\psi)}\biggr) \\
 =& \frac{ {t}^* \pd {u_X}_* \Bigl(\sum_{\nu} {\varphi_{\nu}}_* \frac{1}{s!} 
 \bigl(e_{\nu} \cup u_{\nu}^* \prod_j \prod_{k_j=0}^{d_0 l_j} (l_j H +k \hbar)
 \bigr) \cap [\gbar_{0,0}]^{\vir} \Bigr) } 
  {\prod_{k=1}^{d} (h+k \hbar)^{r+1}} \\
 =&\frac{ {t}^* \sum_{\nu} \frac{1}{s!} {\psi_{\nu}}_* \bigl(
  {u_{\nu}^{\vir}}_* e_{\nu} \cup \prod_j \prod_{k_j=0}^{d_0 l_j}
  (l_j H +k \hbar) \bigr)} {\prod_{k=1}^{d} (h+k \hbar)^{r+1}}. \\
\end{split}
\end{equation}
Here we have used \eqref{e:vl3} and the left square of the above commutative 
diagram \eqref{e:44}. 

{\small
Similar results holds in concave case:\footnote{We have chosen to use 
small font for the discussion in the convex case in this subsection.} 
\begin{equation*} 
 \begin{split} 
 &i_* J_X^E(\beta,\nu)\\
 = &{\ev^{\vir}_P}_* \left( \frac{j_X^* \sum_{\nu} {\varphi_{\nu}}_* 
  \Bigl(e_{\nu} \cup 
   u_{\nu}^* \prod_j \prod_{m=1}^s \xi_m \prod_{k_j=d_0 l_j+1}^{-1} 
  (l_j H +k \hbar) \cap [D_{\nu}]^{\vir} \Bigr)} {\hbar (\hbar-\psi)}\right) \\
 =&\frac{ {t}^* \sum_{\nu}\frac{1}{s!}{\psi_{\nu}}_* \biggl({u_{\nu}^{\vir}}_* 
   e_{\nu} \cup \prod_j \Bigl( \xi_1 \cdots \xi_s \prod_{k_j=d_0 l_j}^{-1}
  (lH +k \hbar) \Bigr)\biggr)} 
  {\prod_{k=1}^{d} (h+k \hbar)^{r+1}}.
 \end{split}
\end{equation*}
}

Apply the correspondence of residues \eqref{e:vl2} to the right
commuting triangle of \eqref{e:43}
\[
 \xymatrix{(\pp^1)^s \times \pp^r_{d_0} \ar[rr]^{\psi_{\nu}} &&\pp^r_d \\
    \pp^r \ar[u]^{t_{\nu}} \ar[rr]^{\operatorname{id}} &&\pp^r \ar[u]_{t} }
\]
one has
\begin{equation} \label{e:46}
\begin{split}
 &\frac{t^* {\psi_{\nu}}_* w}{\prod_{k=1}^{d} (h +k \hbar)^{r+1}} \\
=&\frac{t_{\nu}^* w} {\hbar^s \prod_{k=1}^{d_0}(h+k \hbar)^{r+1}},
\end{split}
\end{equation} 
Then equation \eqref{e:44} becomes, by \eqref{e:46}, 
\begin{equation} \label{e:47}
 i_* J_X^E(\beta,\nu)= \sum_{\nu} \frac{ \prod_j \prod_{k=0}^{d_0 l_j} 
  (l_j h +k \hbar) \cup t_{\nu}^* {u^{\vir}_{\nu}}_* e_{\nu}}
  { s! \hbar^s \prod_{k=1}^{d_0} (h+k \hbar)^{r+1}}.
\end{equation}
{\small
Similarly in the concave case ($l_j <0$)
\begin{equation} \label{e:48}
 i_* J_X^E(\beta,\nu)= \sum_{\nu} \frac{\prod_j \xi_1 \cdots \xi_s 
  \prod_{k=d_0 l_j+1}^{-1} (l_j h +k \hbar) 
  \cup t_{\nu}^* {u^{\vir}_{\nu}}_* e_{\nu}}
	{ s! \hbar^s \prod_{k=1}^{d_0} (h+k \hbar)^{r+1}}.
\end{equation}
}
\emph{Therefore it remains to compute ${u_{\nu}^{\vir}}_* e_{\nu}$.}

First let us work on the hypersurface case, i.e.~$E=L$ and $L$ convex. 
Define $L_{\beta}^k:= i_G^* (\E_d^k)$ for the notational convenience.
From Proposition~\ref{p:0}
\begin{equation*} 
  e_{\nu}= \prod_{m=1}^s \sum_{k=1}^{d_m l} (-k_m) p_{\beta_m}^* 
  \e (L_{\beta_m}^{k_{m}+1}) \cup u_{\nu}^* \e (\F^{k_m}_{\Sigma_m}/
	\F^1_{\Sigma_m})
\end{equation*}
where $\Sigma_m = d_0 +d_1 +\cdots+d_{m-1}$. Again, to simplify
the notations, let us start with the divisorial strata, 
i.e.~$\nu=(\beta_0,\beta_1)$. In this case,
\begin{equation*} 
 e_{\nu}= \sum_{k=1}^{d_1 l} -(k) p_{\beta_1}^* 
  \e (L_{\beta_1}^{k+1}) \cup u_{\nu}^* \e (\F^{k}_{\beta_0}/
	\F^1_{\beta_0}).
\end{equation*}
Therefore
\begin{equation*}
  {u^{\vir}_{\nu}}_* e_{\nu}= \sum_{k=1}^{d_1 l} -(k) {u_{\nu}^{\vir}}_* 
 p_{\beta_1}^* \e (L_{\beta_1}^{k+1}) \cup  \e (\F^{k}_{\beta_0}/
  \F^1_{\beta_0}).
\end{equation*}
It is easy to see that
\begin{equation*} 
 \e (\F^k_{\beta_0}/\F^1_{\beta_0})= \prod_{m=1}^{k-1} (l(h +d_0 \hbar) +
  m \xi)
\end{equation*}
where $\xi$ is the equivariant point class of $\pp^1$. Thus we only have
to know ${u^{\vir}_{\nu}}_* p_1^* \e (L_{\beta_1}^{k+1})$, which is in turn
$\pd {u_{X,1}}_* {p_0}_* (p_1^* \e (L_{\beta_1}^{k+1})\cap 
[{D}_{\nu}]^{\vir})$ because the morphism $u_{\nu}: D_{\nu} \to 
\pp^1 \times \pp^r_d$ factors through $p_0$
\[ 
   D_{\nu} \overset{p_0}{\longrightarrow} G_{0,1}(X,\beta) 
  \overset{u_X,1}{\longrightarrow} \pp^1 \times \pp^r_d.
\]  
For this, consider the fibre square:
\begin{equation} \label{e:53}
\xymatrix{&{D}_{\nu} \ar[ld]_{p_0} \ar[rd]^{p_1} \\
 \gbar_{0,1}(X,\beta_0) \ar[rd]_{\ev_G} 
  &&\mbar_{0,1}(X,\beta_1) \ar[ld]^{\ev_M}\\
  &X}
\end{equation}
so that 
\[ 
 {p_0}_* \bigl( p_1^* \e (L_{\beta_1}^{k+1}) \cap 
 [{D}_{\nu}]^{\vir} \bigr)= \bigl(\ev_G^* {\ev_M^{\vir}}_* 
 \e (L_{\beta_1}^{k+1}) \bigr) \cap [\gbar_{0,1}(X,\beta_0)]^{\vir}. 
\]
However, by dimension counting, 
${\ev_M^{\vir}}_* \e (L_{\beta_1}^{k+1})$ has the cohomological degree 
$2-k+\la c_1(L)-c_1(X), \beta_1 \ra$, which is at most
one for $k=1$ and at most zero for $k=2$. It vanishes otherwise.
It is also obvious that if $\la c_1(L)-c_1(X), \beta \ra \le -2$ for any
$\beta$, i.e.~the case Fano of index $\ge 2$, then all boundary contributions
vanish. 

From our assumption that $i:X \to \pp^r$ induces an isomorphism 
$i^*: NS(\pp^r) \to NS(X)$, the above degree one algebraic cohomology
classes lies in the image of $i^*$ and we may conclude that
\[ \ev_{G}^* {\ev_M^{\vir}}_* \e (L_{\beta_1}^{k+1})
 = i_G^* \ev_{\gbar_{0,1}(\pp^r,d_0)}^* (c^0_{\nu}+ c^1_{\nu} h)  
\] 
for $c^i \in \qq$.
However, as shown in Lemma~\ref{l:1} that 
\begin{equation} \label{e:53-1}
  \ev_{\gbar_{0,1}(\pp^r,d_0)}^* h= u^* (H + d_0 \xi) - \epsilon, 
\end{equation}
where $\epsilon$ is the exceptional divisor of $u$. Summing up, the boundary 
term for the strata $\nu=(\beta_0,\beta_1)$ is  (see \eqref{e:47})
\begin{equation*}
\begin{split} 
 &i_* J^L_X(\beta,\nu) 
  = \frac{\prod_{k=0}^{d_0 l} (l h +k \hbar)}{\prod_{k =1}^{d} 
 (h +k \hbar)} t_{\nu}^* {u_{\nu}^{\vir}}_* e_{\nu} \\
 = &i_* \bigl( J_{X}(\beta_0) \cup H^L_{\beta_0} \bigr) \cup 
  \prod_{k=0}^{d_0 l}   
 \bigl(c_{\nu}^0(l(h+d_0 \hbar) +\hbar) + c_{\nu} (h_i+d^i_0 \hbar) \bigr). 
\end{split}
\end{equation*}
Since $c^0,c^i$ are independent of $\beta_0$ from the above discussion,
this proved the part \emph{(a)} of the following proposition.

\begin{proposition} \label{p:1}
Let $L$ be a convex line bundle on $X$ induced from $P$, and let 
$\lambda_{\nu} (p,\hbar)$ be defined by the following  formula
\begin{equation} \label{e:55}
 i_* \bigl( J_{X}(\beta_0) \cup H^L_{\beta_0} \cup \lambda_{\nu}(p,\hbar)\bigr)
 := \frac{\prod_{k=0}^{d_0 l} (l h+k \hbar) \cup t_{\nu}^* 
 {u^{\vir}_{\nu}}_* e_{\nu}} {\prod_{k=1}^{d_0} (h+k \hbar)^{r +1}} ,
\end{equation} 
which is $i_* J_X^L(\beta, \nu)$.

(a) If $\nu=(\beta_0,\beta_1)$, then $\lambda_{\nu}(p,\hbar)$ is linear 
and satisfies 
\[ \lambda_{(\beta_0,\beta_1)}(p,\hbar)
 =\lambda_{(0,\beta_1)}(p+ d_0  \hbar, \hbar).
\]

(b)
for $\nu=(\beta_0,\beta_1, \cdots,\beta_s)$, we have
\begin{equation} \label{e:56}
\lambda_{\nu}(p,\hbar)=\lambda_{(\Sigma_1,\beta_1)}(p,\hbar) 
\lambda_{(\Sigma_2,\beta_2)}(p,\hbar) \cdots 
\lambda_{(\Sigma_{s},\beta_s)}(p,\hbar),
\end{equation}
where $\Sigma_m = \sum_{a<m} \beta_a$.
\end{proposition}

\begin{proof} (of part \emph{(b)})
The equation \eqref{e:56} requires only to compute $t_{\nu}^* {u_{\nu}}_*
e_{\nu}$, with $e_{\nu}$ described by \eqref{e:47}
\begin{equation*} 
  {u_{\nu}}_* e_{\nu}
 = \prod_m \sum_{k_m} (-k_m) {u_{\nu}}_* p_{\beta_m}^* 
 \e (L_{\beta_m}^{k_m+1}) \cup \e (\F^{k_m}_{\Sigma_m}/\F^1_{\Sigma_m}).
\end{equation*}
Using the same argument 
to the following diagram
\begin{equation} \label{e:59}
\xymatrix{&{D}_{\nu} \ar[ld]_{p_0} \ar[rd]^{p_{M}} \\
  \gbar_{0,s}(X,\beta_0) \ar[rd]_{\ev_G} &&\prod \mbar_{0,1}(X,\beta_m) 
 \ar[ld]^{\ev_{M}}\\
  &(X)^s}
\end{equation}
one obtains that 
\begin{equation*} \label{e:60}
 \e (\F^{k_m}_{\Sigma_m}/\F^1_{\Sigma_m})= \prod_{a=1}^{k_m-1} \biggl(lh+
 l \Sigma_m \hbar + \sum_{b=1}^{m-1} l d_{b} \xi_{b} + a \xi_m \biggr).
\end{equation*}
Dimension counting and \eqref{e:53-1} gives us the proof of part \emph{(b)}.
\end{proof}

{\small
In the concave case a similar modification goes through.

\begin{proposition} \label{p:2}
Let $L$ be a concave line bundle and let $\lambda_{\nu}$ be defined as in
\eqref{e:55}. Then part (a) and (b) hold. Furthermore, $\lambda_{\nu}$ depends
only on $\hbar$.
\end{proposition}

\begin{proof}
From \eqref{e:48} for $J^L_P(\beta,\nu)$, we only have to compute 
${u_{\nu}^{\vir}}_* e_{\nu}$. A straightforward modification, by fibre
square \eqref{e:53} and dimension counting, will lead to the conclusion that 
only degree zero (constant) terms survives in ${\ev_M^{\vir}}_* 
\e (L_{\beta_1}^{k+1})$. The difference is that
the cohomological degree of ${\ev_M^{\vir}}_* \e (L_{\beta_1}^{k+1})$
is only $-k+\la -c_1(L)-c_1(X), \beta_1 \ra$, because the rank of $L^{k+1}$
in this case is one less than that of $L^{k+1}$ in convex case.

Therefore $\lambda_{\nu}= constant \cdot t_{\nu}^* \prod_{m=1}^s
\xi_m = constant \cdot \hbar^s$. This allows us to set 
$\lambda_{(\beta_0,\beta_1)}= constant \cdot \hbar$.
\end{proof}
}

\begin{corollary} \label{c:4}
The same is true when the rank of $E$ is greater than 1.
\end{corollary}

The proof uses the inductive property of $e_{\nu}$ and is completely
analogous to the above argument. 

Summarizing the above discussion, we have
\begin{theorem} \label{t:5}
\begin{equation} \label{e:61}
 i_* J_X^E = i_* \sum_{\beta} q^{\beta} \sum_{\nu} \frac{
 J_X({\beta}_0) \cup H^E_{\beta_0} \cup \prod_{m=1}^s \lambda_{{\beta}_m} 
 (p + \Sigma_m \hbar, \hbar)}{s! \hbar^s}
\end{equation}
such that $\lambda_{{\beta}_m}(p,\hbar)$ are linear form in $p,\hbar$ and
is independent of the total degree $\beta$ in the background.
\end{theorem}
\begin{proof}
Set $\lambda_{\beta_m}(p,\hbar)= \lambda_{(0,\beta_m)}(p,\hbar)$ and apply
Propositions~\ref{p:1} and \ref{p:2}.
\end{proof}

\subsection{Mirror transformation} \label{ss:4.4}
\emph{In this subsection we will omit the push-forward symbol $i_*$.
All equalities are assumed to hold after pushing-forward to $\pp^r$.}

Now we are ready to find a change of variables (mirror transformation)
between $J_X^E$ and $I_X^E$. Note that we will use $\beta$ both
as an element in $H_2(X)$ and a number. 
In the case $Pic(X)=1$, this should not cause confusion.
\begin{equation*} 
\begin{split}
  &t_0 \mapsto t_0 + f_0(q) \hbar  , \\
  &t \mapsto t + f_1(q), 
\end{split}
\end{equation*}
where $f_0(q), f_1(q)$ are formal power series in $q$. 
After the change of variables, $e^{(t_0 + p t)/\hbar}I_X^E(q, \hbar)$ 
becomes
\begin{equation} \label{e:63}
  e^{f_0(q)+\frac{1}{\hbar} p f_1(q)} e^{\frac{1}{\hbar} (t_0 + p t)} 
  I_X^E(q e^{f_1(q)}, \hbar),
\end{equation}
\emph{which we would like to equate to $J_X^E(q,\hbar)$}.

Let us set $f_1(q)= \sum_{\beta} a_{\beta} q^{\beta}$ and $f_0(q)= 
\sum_{\beta} b_{\beta} q^{\beta}$. Then \eqref{e:63} can be expanded as:
\begin{equation} \label{e:72}
\begin{split}
 &e^{f_0(q) +\frac{1}{\hbar} h f_1(q)} e^{\frac{1}{\hbar}
   (t_0 +  p t)} I_X^E(q e^{f_1(q)}, \hbar) \\
 =&\sum_{\beta_0} q^{\beta_0} e^{f_0(q)+
  (\frac{p}{\hbar}+{\beta}_0) f_1(q)} I_{X}^E ({\beta}_0) \\
 =&\sum_{\beta} q^{\beta} \sum_{\sum_{m=0}^s {\beta}_m ={\beta}} \frac{1}{s!} 
   I^E_{X}({\beta}_0) \prod_{m=1}^s \Bigl( b_{{\beta}_m} 
  + a_{{\beta}_m} ({\beta}_0 + \frac{p}{\hbar}) \Bigr).
\end{split}
\end{equation}

In order to prove Theorem~\ref{t:1} we need to find $a_{\beta},b_{\beta}$ so
that \eqref{e:72} is equal to $J_X^E(q,\hbar)$. The following simple lemma 
in \cite{AB} is useful:

\begin{lemma} \label{l:2}
Let
\begin{equation} \label{e:73}
 Q(q)= \sum_{{\beta}} q^{\beta} \sum_{\sum_{m=1}^s {\beta}_m ={\beta}, 
  d_m \ne 0} \frac{1}{s!} 
 \prod_{m=0}^s \Bigl(y_{{\beta}_m}+ x_{{\beta}_m}{B_{m-1}} \Bigr),
\end{equation}
where $B_{m-1}:= \sum_{k=1}^{m-1} \beta_k$, $(B_0=0)$. 
Then $\log (Q(q))$ is a linear function of $y_{\beta}$.
\end{lemma}


\begin{corollary} \label{c:5}
Use the notation from the above lemma. Define $z_{\beta}$ by 
\begin{equation} \label{e:74}
Q(q)= \exp (\sum_{\beta \ne 0} z_{\beta} q^{\beta}). 
\end{equation}
Then
\begin{equation} \label{e:75}
 z_{\beta} = \sum_{\sum_{m=1}^s {\beta}_m={\beta}} \frac{1}{s!} y_{{\beta}_1} 
  \prod_{m=1}^s x_{{\beta}_m} B_{m-1}
\end{equation}
\end{corollary}

\begin{proof}
Expand equations \eqref{e:74} and \eqref{e:73}. The $q^{\beta}$ term is
\begin{equation} \label{e:76}
\sum_{\sum_{m=1}^s \beta_m =\beta} \frac{1}{s!} \prod_{m=1}^s 
z_{{\beta}_m} = \sum_{\sum_{m=1}^s \beta_m = \beta} \frac{1}{s!} \prod_{m=1}^s
(y_{{\beta}_m} + x_{{\beta}_m} B_{m-1}). 
\end{equation}
By Lemma~\ref{l:2} the right hand side of \eqref{e:76} is linear with respect 
to $y_{\beta}$. The linear in $y_{\beta}$ term on the RHS is exactly 
\eqref{e:75}.
\end{proof}

Back to the proof of Theorem~\ref{t:1}. We wish to equate \eqref{e:61} to 
\eqref{e:72}. Set $\lambda_{\beta} = c^0_{\beta} \hbar + c^1_{\beta} p$, 
then \eqref{e:61} can be written as 
\begin{equation*} 
\begin{split}
 J_X^E(q) = &\sum_{{\beta}_0} \sum_{\sum_{m=1}^s {\beta}'_m={\beta}'} 
 q^{{\beta}'+{\beta}_0} J_{X}({\beta}_0)
 \Bigl( \prod_j H_{\beta_0}^{L_j} \Bigr)  \\
  &\qquad \qquad \frac{1}{s!} \prod_{m=1}^s 
  \biggl( \Bigl(c^0_{{\beta}'_m} + c^1_{{\beta}'_m} (\frac{p}{\hbar} 
  +{\beta}_0) \Bigr) + c^1_{{\beta}'_m} {B_{m-1}} \biggr) \\
 =& \sum_{\beta_0} q^{\beta_0} J_X(\beta_0) 
  \Bigl(\prod_j H^{L_j}_{\beta_0} \Bigr) \\
 &\qquad \qquad \sum_{\sum_{m=1}^s \beta'_m = \beta'} q^{\beta'}\frac{1}{s!}
  \prod_{m=1}^s \biggl(\Bigl(c^0_{{\beta}'_m}+c^1_{{\beta}'_m}(\frac{p}{\hbar} 
  +{\beta}_0) \Bigr) + c^1_{{\beta}'_m} B_{m-1} \biggr).
\end{split}
\end{equation*}
Apply Corollary~\ref{c:5}, we find a new variables $z_{{\beta}'}$ such that
$z_{{\beta}'}$ are linear in $c^0_{{\beta}'_m} + c^1_{{\beta}'_m} 
(\frac{p}{\hbar} +{\beta}_0)$ and polynomial in $c^1_{{\beta}'_j}$. 
Therefore there are constants $b_{{\beta}'}, a_{{\beta}'}$ such that 
$z_{{\beta}'}= b_{{\beta}'} + a_{{\beta}'} (\frac{p}{\hbar}+{\beta}_0)$, 
and $b_{{\beta}'}, a_{{\beta}'}$ are independent of ${\beta}_0$. 
This implies that
\begin{equation*} 
 J_X^E(q) = \sum_{\sum_{m=0}^s {\beta}_m={\beta}} q^{\beta} I_{X}^E({\beta}_0) 
 \frac{1}{s!} \prod_{m=1}^s \Bigl(b_{{\beta}_m} + a_{{\beta}_m} 
 (\frac{p}{\hbar} +{{\beta}_0}) \Bigr),
\end{equation*}
which is exactly \eqref{e:72}. One can easily see from the above proof that
the case $E$ being a direct sum of convex and concave bundles requires little
modification. Our proof of the Theorem~\ref{t:1} is therefore complete.

\subsection{Proof of Theorem~\ref{t:2}}
The proof of case~\ref{convex} is very easy. As we have seen in 
\S~\ref{ss:boundary} that the use of fibre square \eqref{e:53} \eqref{e:59}
and dimension counting guarantees that ${u_{\nu}^{\vir}}_* e_{\nu} =0$.
Therefore $J_X^E(\beta,\nu)=0$.

The proof of case~\ref{concave} goes a slightly different way.\footnote{It is
also possible to prove the convex case in this way. \label{f:9}}.
As we have seen in the proof of Proposition~\ref{p:2} that ${u_{\nu}^{\vir}}_* 
e_{\nu}$ contains only constant terms. Therefore the boundary contribution
should come from
\begin{equation*} 
 \gbar_{0,s}(X,\beta_0) \overset{u_s}{\to} (\pp^1)^s \times P_{d_0} 
 \overset{\psi_{\nu}}{\to} \pp^r_d,
\end{equation*}
and should be of the form
\begin{equation*} 
 {\psi_{\nu}}_* \left( (constant) ({u_s}_* 1) \cup \prod_{j=1}^{\rk(E)} 
 \xi_1 \cdots \xi_s \prod_{k_j=\la c_1(L_j),\beta_0 \ra +1}^{-1} 
  (l_j H +k_j \hbar) \right),
\end{equation*}
which has cohomological dimension greater than $J_X^E(\beta,\text{main})$.
Therefore $constant=0$.

The proof of direct sum case is the combination of the above arguments.

\subsection{Proof of Corollary~\ref{c:qsd}}
To simplify our notations, set 
\[
 J':=i_* \e(E) e^{\frac{1}{\hbar} (t_0 + p t)} J^E_X(q,\hbar),
\]
and 
\[
 J'_{\vee} :=i_* \frac{1}{\e(E)} 
  e^{\frac{1}{\hbar}(t_0 + p t)} J^{E^{\vee}}_X(q,\hbar).
\]
Notice that $J'_{\vee}$ is actually well defined as  $J^{E^{\vee}}_X(q,\hbar)$
always has a factor $\e(E^{\vee})$. More precisely, consider the exact sequence
\[
 0 \to R^0 \pi_* (\ev^*(E^{\vee}) \otimes \so(-x_1)) \to R^0 \pi_* 
\ev^*(E^{\vee}) \to \ev^* E^{\vee} \to 0.
\]
Then in fact
\[
\begin{split}
J'_{\vee} =& i_* e^{\frac{1}{\hbar}(t_0 + p t)}
 \sum_{\beta} q^{\beta} \ev^{\vir}_* \left(
\frac{\e\Bigl(R^0 \pi_* \bigl(\ev^*(E^{\vee}) 
 \otimes \so(-x_1) \bigr) \Bigr)}{\hbar(\hbar -\psi)} \right) \\
&=  i_* e^{\frac{1}{\hbar}(t_0 + p t)} 
  J^{E^{\vee}}_X \frac{1}{\e(E^{\vee})}.
\end{split}
\] 
Similar interpretation holds for $J'$: it is the $J$-function of the
bundle $R^1 \pi_* (\ev^*(E) \otimes \so(-x_1))$, multiplied by 
$e^{\frac{1}{\hbar}(t_0 + p t)}$. Notice that the ranks of the two
bundles $R^0 \pi_* (\ev^*(E^{\vee}) \otimes \so(-x_1))$ and 
$R^1 \pi_* (\ev^*(E) \otimes \so(-x_1))$ are the same.

It follows from Theorems~\ref{t:1} and \ref{t:2}
\[
  J'= I' := e^{\frac{1}{\hbar}(t_0+ p t)} 
  \sum_{\beta} q^{\beta} J_X(\beta) \cup \prod_j 
  \prod_{k_j=1}^{\la c_1(L_j), \beta\ra} (c_1(L) +k \hbar)
\]
and 
\[
 J'_{\vee} \sim I'_{\vee}:= (-1)^{\rk(E)} e^{\frac{1}{\hbar}(t_0+ p t)} 
  \sum_{\beta} q^{\beta}  J_X(\beta) \cup \prod_j 
 \prod_{k_j= \la c_1(L_j), \beta\ra +1}^{0} (c_1(L) +k \hbar).
\]
Notice the change of limits in the products of $k$ due to $\e(E)$ and 
$\e(E^{\vee})$ and the possible sign 
coming from the ratio of $\e(E)$ and $\e(E^{\vee})$.

The difference between $I'$ and $I'_{\vee}$ can be easily figured out by using
the quantum differential equation \cite{ABG1}.

\begin{claim}
There is a power series $\phi(q)$ such that 
\begin{equation} \label{e:qsd1}
\phi(q) I'= I'_{\vee}.
\end{equation}
\end{claim}

We now sketch the proof of \eqref{e:qsd1}. 
Fix a basis $\{e_a\}$ of the vector space $H^*(X)$ with $e_1=1$. Let 
\[ 
 S_{ab}:= \sum_{\beta} q^{\beta} \bigg(e_a, e^{\frac{1}{\hbar}(t_0 + p t)}
 \frac{e_b}{\hbar -\psi} \bigg)_{0,2,\beta}
\]
be a matrix of genus zero two-point Gromov--Witten invariants via the 
``fundamental class'' 
$\e(E_{\beta}) \cup \e(E) \cap [\mbar_{0,2}(X,\beta)]^{\vir}$.
Note that we have identified $t_i=\log q_i$. Givental's theory of quantum 
differential equation says 
\begin{equation} \label{e:qde}
 \hbar \frac{\p}{\p t_i} S= p_i \circ S, 
\end{equation}
where $p_i \circ$ is the quantum multiplication matrix.
By the (virtual) dimension counting $J_X$ is a polynomial in $\hbar^{-1}$ of 
the form $1+O(\hbar^{-\la c_1(E), \beta \ra})$. Therefore 
\[
 I'=J'=(1+O(\hbar^{- \rk(E)})) e^{\frac{1}{\hbar}(t_0 + p t)}.
\]
Notice that $S_{1a}=\la J', e_a \ra$. This implies that the first row
$S_{1a}$ of $S$ has the same order in $\hbar^{-1}$. 

Set $v_j:=c_1(L_j)$ and $\p_j$ be the directional derivative on the direction 
$v_j$. First of all $(\prod_{j=1}^{\rk(E)} \p_j) ({I'}/{\e(E)}) = I'_{\vee}$ 
by simple derivations. Now apply induction on the order $m$ of the differential
operator $\mathcal{D}_m:=\prod_{j=1}^m \hbar \p_j$. I claim that
\begin{equation} \label{e:qsd2}
  \mathcal{D}_m(S)= \left(\prod_{j=1}^m v_j \right) \circ S
\end{equation}
for $m < \rk(E)$.

The case $m=1$ is true by \eqref{e:qde}.
Suppose that \eqref{e:qsd2} is true 
for some $m < \rk(E)-1$. This implies that 
$\mathcal{D}_m S= (\prod_j^m v) \circ + O(\hbar^{-1})$. Therefore 
$((\prod_j^m v)\circ)_{1a}$ are independent of $t$ from LHS. Now differentiate
one more time:
\begin{equation} \label{e:qsd3}
 \mathcal{D}_{m+1} S= \left(\hbar \p_{m+1} \Bigl(\prod_{j=1}^m v_j \Bigr) 
  \circ \right)S +\left( \Bigl(\prod_{j=1}^m  v_j\Bigr) \circ 
	v_{m+1} \circ \right) S.
\end{equation}
and the first row of RHS will be $\la (\prod_{j=1}^m  v_j) \circ v_{m+1},
e_a \ra + O(\hbar^{-1})$. Because the first row of the LHS modulo 
$\hbar^{-1}$ is independent of $t$, we have again 
$(\prod_{j=1}^{m}  v_j) \circ v_{m+1}= \prod_{j=1}^{m+1} v_j$. Thus 
\eqref{e:qsd2} holds. If $m=\rk(E)-1$, then the first row of LHS of 
\eqref{e:qsd3} would be certain power series $\phi(q)$ of $q$ 
(mod $\hbar^{-1}$). Since the first row of first term on RHS still vanishes 
(mod $\hbar^{-1}$), we have the first row of \eqref{e:qsd3} equal to
\[
\begin{split}
 \phi(q) I'
=&\phi(q) \left(\prod_{j=1}^{\rk(E)}v_j \right) \left(\frac{I'}{\e(E)}\right)\\
 =&\left(\prod_{j=1}^{\rk(E)} (v_j \circ) \right) I' =\mathcal{D}_{\rk(E)} I'\\
 = &\pm I'_{\vee}.
\end{split}
\]

\end{document}